\documentclass[12pt]{amsart}
\usepackage{graphicx}
\newtheorem{thm}{Theorem}[section]

\newtheorem{prop}[thm]{Proposition}
\theoremstyle{definition}
\newtheorem{defn}[thm]{Definition}
\theoremstyle{remark}

\theoremstyle{definition}
\newtheorem{alg}[thm]{Scheme}
\numberwithin{equation}{section}

\begin{document}

\title[Closed queueing networks in semi-Markov environment]
{Large closed queueing networks in semi-Markov environment
and their application}%
\author{Vyacheslav M. Abramov}%
\address{School of Mathematical Sciences, Monash University, Building 28M, Wellington road, Clayton, VIC 3800, Australia}%
\email{vyacheslav.abramov@sci.monash.edu.au}%

\subjclass{60K25, 60K30, 60H30, 60H35}%
\keywords{Closed queueing network; Random environment; Martingales
and Semimartingales; Skorokhod reflection principle}%

\begin{abstract}
The paper studies closed queueing networks containing a server
station and $k$ client stations. The server station is an infinite
server queueing system, and client stations are single-server
queueing systems with autonomous service, i.e. every client
station serves customers (units) only at random instants generated
by a strictly stationary and ergodic sequence of random variables.
The total number of units in the network is $N$. The expected
times between departures in client stations are $(N\mu_j)^{-1}$.
After a service completion in the server station, a unit is
transmitted to the $j$th client station with probability $p_{j}$
$(j=1,2,\ldots,k)$, and being processed in the $j$th client
station, the unit returns to the server station. The network is
assumed to be in a semi-Markov environment. A semi-Markov
environment is defined by a finite or countable infinite Markov
chain and by sequences of independent and identically distributed
random variables. Then the routing probabilities $p_{j}$
$(j=1,2,\ldots,k)$ and transmission rates (which are expressed via
parameters of the network) depend on a Markov state of the
environment. The paper studies the queue-length processes in
client stations of this network and is aimed to the analysis of
performance measures associated with this network. The questions
risen in this paper have immediate relation to quality control of
complex telecommunication networks, and the obtained results are
expected to lead to the solutions to many practical problems of
this area of research.

\end{abstract}
\maketitle

\tableofcontents
\newpage
\section{Introduction}
We consider closed queueing networks containing a server station
and $k$ client stations. The server station is an infinite server
queueing system with identical servers. Client stations are
single-server queueing systems with an autonomous service
mechanism, where customers (units) are served only at random
instants generated by a strictly stationary and ergodic sequence
of random variables.

Queueing systems with an autonomous service mechanism were
introduced and originally studied by Borovkov \cite{Borovkov
1976}, \cite{Borovkov 1984}. The formal definition of these
systems in the simplest case of single arrivals and departures is
as follows. Let $A(t)$ denote an arrival point process, let $S(t)$
denote a departure point process, and let $Q(t)$ be a queue-length
process, and all of these processes are started at zero
($A(0)=S(0)=Q(0)=0$). Then the autonomous service mechanism is
defined by the equation:
\begin{equation*}
Q(t)=A(t)-\int_0^\infty \mathbb{I}\{Q(s-)>0\}\mbox{d}S(s).
\end{equation*}
The queueing systems with an autonomous service mechanism have
been studied in many papers (e.g. Abramov \cite{Abramov 2000},
\cite{Abramov 2004}, \cite{Abramov 2006}, Fricker \cite{Fricker
1986}, \cite{Fricker 1987}, Gelenbe and Iasnogorodski
\cite{Gelenbe and Iasnogorodski 1979}). The structure of queueing
systems or networks with autonomous service and their analysis is
much easier than those structure and analysis of ``usual" systems
with generally distributed service times. Queueing systems and
networks with autonomous service, because of their simple
construction, are studied under general settings on dependent
inter-arrival and inter-departure times, and their analysis is
often based on methods of stochastic calculus and the theory of
martingales. The corresponding results for usual Markovian
queueing systems or networks follow as a particular case of the
corresponding results for queueing systems or networks with
Poisson input and autonomous service. For different applications
of queueing systems (networks) with autonomous service see e.g.
\cite{Abramov 2000}, \cite{Abramov 2004}, \cite{Abramov 2005} and
\cite{Abramov 2006}.

The assumption that the queueing mechanism is autonomous
substantially simplify the analysis. However, according to
sample-path results of \cite{Abramov 2005} this assumption can be
removed. So, all of the main results remain valid for quite
general client/server networks without the special assumption that
the service mechanism is autonomous.

In the present paper we study client/server networks in a
semi-Markov environment. There has been an increasing attention to
queueing systems in a random environment in the literature (e.g.
\cite{Baccelli and Makovsky 1986}, \cite{Boxma and Kurkova 2000},
\cite{D'Auria 2007}, \cite{Helm and Waldmann 1984}, \cite{Krieger
et al 2005}, \cite{O'Cinneide and Purdue 1986}, \cite{Purdue 1974}
and others). However, most of these papers mainly develop the
theory and remain far from real-world applications.

The aim of the present paper is twofold. First, we establish new
theoretical results for client/server networks in semi-Markov
environment describing the behavior of queue-length processes of
this network. Second, we show how these theoretical results can be
applied to solve real-world problems. Some of these problems are
solved in the present paper. Other ones will be solved in the
future.

The model of the network, which is considered in this paper, is
very close to the models considered in \cite{Kogan and Liptser
1993} and \cite{Abramov 2000} (see Figure 1). For a more general
construction of network with two types of node and multiple
customer classes see \cite{Abramov 2004}.

The description of the present model is based on that of the model
of \cite{Abramov 2000}. For other papers studying the models of
client/server networks see also \cite{Berger Bregman and Kogan
1999}, \cite{Kogan 1992}, \cite{Kogan Liptser and Smorodinskii
1986}, \cite{Krichagina Liptser and Puhalskii 1988},
\cite{Krichagina and Puhalskii 1997}, \cite{McKenna and Mitra
1982}, \cite{Pittel 1979}, \cite{Whitt 1984} and other papers.

\begin{figure}
\includegraphics[width=15cm,height=20cm]{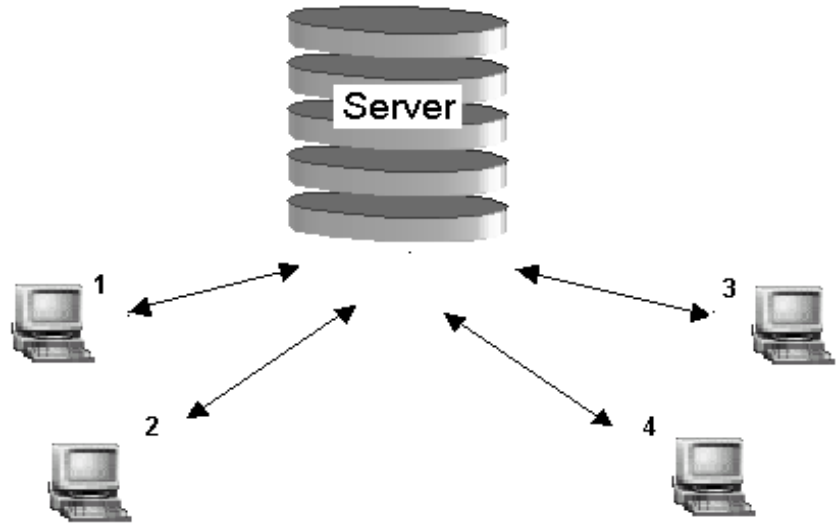}
\caption{An example of client/server network topology}
\end{figure}

The previous assumptions of \cite{Abramov 2000}, that are also
used in the present paper, are repeated below.

The departure instants in the $j$th client station
($j=1,2,\ldots,k$) are denoted by $\xi_{j,N,1}$,
$\xi_{j,N,1}+\xi_{j,N,2}$,
$\xi_{j,N,1}+\xi_{j,N,2}+\xi_{j,N,3}$,\ldots, and

\smallskip

$\bullet$ each sequence $\{\xi_{j,N,1}, \xi_{j,N,2},\ldots \}$
forms a strictly stationary and ergodic sequence of random
variables.

\smallskip
The corresponding point process associated with departures from
the client station $j$ is denoted
$$S_{j,N}(t) = \sum_{i=1}^\infty
\mathbb{I}\left\{\sum_{l=1}^i\xi_{j,N,l}\leq t\right\}.$$

The total number of units in the network is $N$. The number $N$ is
a large parameter, and we assume that $N$ increases indefinitely.
This means that we assume that $N$ is a series parameter, and the
series of models (with different $N$) are considered on the same
probability space.

It was assumed in \cite{Abramov 2000} that the service time of
each unit in the server station is an exponentially distributed
random variable with a \textit{given} parameter $\lambda$. In the
present paper, the assumption is another. Under the assumption
that an environment is random  (Markov or semi-Markov), the
parameter $\lambda$ is not longer a \textit{constant value}. It is
a random variable, taking values in dependence of environment
states. (By environment states we mean the states of the
corresponding Markov chain, which is described later.) The same is
noted regarding the routing probability matrix.

There were the standard assumptions in \cite{Abramov 2000} that
after a service completion at the server station, a unit was
transmitted to the client station $j$ with probability $p_{j}$,
$p_{j}\geq 0$, and $\sum_{j=1}^{k}p_{j}$=1. These assumptions are
not longer valid in the case of the system considered in this
paper. The routing probabilities are assumed to be random, taking
the values in dependence of the environment states as well.

Let us first describe \textit{Markov states} of semi-Markov
environment, and then latter we describe the notion of semi-Markov
environment as well. The network is assumed to be complete in the
following sense.

The server station and $k$ client stations communicate by links.
The number of links is $k$. Therefore, the Markov chain of the
environment states is defined as follows. The number of possible
environment states is assumed to be finite or countable infinite.
These states are denoted by $\mathcal{E}_{i}$, $i=1,2,\ldots$. For
the Markov or semi-Markov environment considered later we shall
also use the notation $\mathcal{E}(t)$. The meaning of this
notation is the state of Markov chain in time $t$. The initial
state is denoted by $\mathcal{E}_0$ or $\mathcal{E}(0)$.
\smallskip

The completeness of the network means that

\smallskip
$\bullet$ for any $j$=1,2,\ldots, $k$ there exists $i_j\geq1$,
such that $p_j(\mathcal{E}_{i_j})>0$.

\smallskip
This assumption is not used in our proofs explicitly.
Nevertheless, it must be mentioned. If $p_j(\mathcal{E}_{i_j})=0$
for all $i_j\geq1$, then the station $j$ becomes isolated and
therefore not representative. In the above assumption by
$p_j(\mathcal{E}_{i})$ we mean the value of probability $p_j$ when
the network is in state $\mathcal{E}_{i}$. It is assumed
additionally that for any $i=1,2,\ldots$, the sum of probabilities
$\sum_{j=1}^{k}p_j(\mathcal{E}_{i})$ =1. The notation with similar
meaning is used for parameter $\lambda$. Namely,
$\lambda(\mathcal{E}_{i})$ is assumed to be strictly positive for
any state $\mathcal{E}_{i}$, however $\lambda_j(\mathcal{E}_{i})$
= $\lambda(\mathcal{E}_{i})p_j(\mathcal{E}_{i})$ can be equal to 0
(because $p_j(\mathcal{E}_{i})$ need not be strictly positive in
general).

\smallskip

The results of \cite{Abramov 2000} are associated with the
asymptotic analysis of a closed client/server network with a
bottleneck station as $N$ increases indefinitely. One of the main
results of \cite{Abramov 2000} was then developed in \cite{Abramov
2004} for networks containing two types of node and multiple
customer classes, where one of client stations was bottleneck. The
results obtained in \cite{Abramov 2000}, \cite{Abramov 2001},
\cite{Abramov 2004} as well as in preceding paper \cite{Kogan and
Liptser 1993} all are a theoretical contribution to the theory of
client/server computer networks with bottlenecks.

The aim of the present paper is another. We follow towards
performance analysis, and are aimed to use the known theoretical
results of \cite{Kogan and Liptser 1993} and \cite{Abramov 2000}.
However, the known theoretical results obtained in these papers
are scanty for their immediate application to real
telecommunication systems, whose parameters can change in time.
Therefore, for the purpose of performance analysis, a substantial
development of the earlier results related to this model is
required. Therefore, before studying the behaviour of queueing
networks in semi-Markov environment, we first study the behaviour
of this network in its special case of piece-wise deterministic
environment. Then the results are extended to the more general
situation of semi-Markov environment.

The paper is organized as follows. In Section \ref{Motivation},
the paper is motivated by formulating the concrete real-world
problems, where the results of the present paper can be applied.
In section \ref{Definitions}, necessary notions of piece-wise
deterministic, Markov and semi-Markov environments as well as
classification of bottleneck stations are defined. In Section
\ref{Semi-Martingale}, the semi-martingale decomposition for the
queue-length process in semi-Markov environment is derived. Then
in Section \ref{Piece-wise deterministic}, the queue-length
processes in piece-wise deterministic environment are studied. In
fact we study ``usual" closed client/server networks under general
assumptions. In Section \ref{Semi-Markov}, the results are
extended to the case of semi-Markov environment. In Section
\ref{Application}, an example of an application of the main
theorems of the present paper to client/server computer network
with failing client stations is considered. In Section
\ref{Discussion}, new problems and monotonicity conditions
associated with these problems for the case of Markov environment
are discussed. In Section \ref{Numerical example}, the example of
numerical study for the simplest network is considered. The
results of the paper are concluded in Section \ref{Concluding
remarks}.

\section{Motivation}\label{Motivation}

In this section we motivate our study by real-world problems
associated with performance analysis of the network, its separate
stations or a subset of those stations. The main results of the
present paper are the subject for many concrete practical problems
having significant value. Two of them are considered in
\cite{Abramov reliability 2007}, \cite{Abramov reliability
2007-2}. Article \cite{Abramov reliability 2007} assumes that
client stations all are identical and subject to breakdowns.
Article \cite{Abramov reliability 2007-2} discusses the similar
problem for not identical client stations. As a client station is
breakdown the parameters of the network are changed. One or other
breakdown leads to bottleneck and risk of a destruction of the
network performance. The aforementioned articles \cite{Abramov
reliability 2007}, \cite{Abramov reliability 2007-2} establish
confidence intervals helping to avoid with a given high
probability that performance destruction. Some results of these
articles are outlined in Section \ref{Application} of this paper.

Another circle of problems is as follows.

 If a computer network operates
for a given fixed time interval (0, $T$), then its performance
characteristic (for example, payment for the increasing a given
level of queue in a given client station, or in a subset of client
stations) depends on the characteristics of this network at the
initial time $t=0$. Using the terminology of the present paper,
these characteristics at the initial time $t=0$ specify the
initial condition of the environment. For one initial condition,
the measure of time that the queue-length is greater than a given
fixed level $L$ is $x$, and we must pay, say $Cx$. For other
initial condition, this measure of time is $y$ and the
corresponding payment is $Cy$.

On the other hand, the cost for initialization the first initial
condition related to the above client station of the network is
$X$, and that cost for initialization the second initial condition
is $Y$.

So, the total expenses in the first case are $X$ + $Cx$, and in
the second case $Y$ + $Cy$. If $X$ + $Cx<Y$ + $Cy$, then we say
that the first \textit{strategy} is more profitable than the
second one. In practical context, the first strategy can mean
\textit{the first type repair} of the network, or a failing client
station, or a subset of failing client stations of the network,
and the second strategy - \textit{the second type repair} of that
network, or a failing client station, or a subset of failing
client stations correspondingly. For example, the first type
repair can contain an additional prophylactical service and
therefore to be more expensive that the second type repair (i.e.
$X>Y$).


\section{Definition of semi-Markov environment and classification
of bottleneck stations}\label{Definitions}

\subsection{Semi-Markov environments}
Let us now define a semi-Markov environment in Mathematical terms.
For each state $\mathcal{E}_i$ of the Markov chain let us define
the sequence of independent and identically distributed random
variables:
\begin{equation}\label{I0}
\zeta_{i,1}, \zeta_{i,2},\ldots
\end{equation}

\begin{defn}\label{Defn-semiMI}Let $\{Z_n\}$ be a Markov chain
with finite or countable infinite states space
$\{\mathcal{E}_i\}$. Let $Z(t)$ be a random process, and let
$\{\sigma_n\}$ be an increasing sequence of random points, i.e.
$0=\sigma_0<\sigma_1<\ldots$. Suppose that the process $Z(t)$ is
defined as follows: $Z(\sigma_n)=Z_n$, and $Z(t)\equiv
Z(\sigma_n)$ for all $\sigma_n\leq t<\sigma_{n+1}$. Suppose also
that the increments $\sigma_{n+1}-\sigma_n$, $n$=0,1,2,\ldots,
coincide in distribution with $\zeta_{i,n}$, where the equality
$\mathcal{E}(\sigma_n)=\mathcal{E}_i$ defines the index $i$. Then
the process $Z(t)$ is called \textit{semi-Markov environment}.
\end{defn}

The above semi-Markov environment $Z(t)$ is assumed to be given on
special probability space $\{\Omega_Z, \mathcal{F}_Z,
\mathbb{P}_Z\}$, which in turn is contained in the common filtered
probability space $\{\Omega$, $\mathcal{F}$,
$\mathbf{F}=(\mathcal{F}_t$), $\mathbb{P}\}$.

The two significant special cases of semi-Markov environment are
as follows.

\begin{defn}\label{Defn-piece-wise-det}
A semi-Markov environment is called piece-wise deterministic, if

\smallskip
$\bullet$ the sequences of \eqref{I0} all are deterministic, i.e.
 $\zeta_{i,v}=z_i$, where $z_i$ is the same constant for all $v$.

\smallskip
$\bullet$ the parameters of networks are not randomly dependent of
the environment states, i.e.
$\lambda_j(\mathcal{E}_0)=\lambda_j^{(0)}$,
$\lambda_j(\mathcal{E}(\sigma_1))=\lambda_j^{(1)}$,\ldots,
$\lambda_j(\mathcal{E}(\sigma_l))=\lambda_j^{(l)}$,\ldots, for all
$j$=1,2,\ldots, $k$, where $\lambda_j^{(0)}$, $\lambda_j^{(1)}$,
\ldots all are non-random constants.

\smallskip

In other words, a piece-wise deterministic environment defines
time dependent closed client/server network with fixed piece-wise
constant parameters depending on time.

\end{defn}

\begin{defn}\label{Defn-MI}
In the case where the sequences \eqref{I0} consist of independent
and exponentially distributed random variables with rates $z_i$,
$i$=1,2,\ldots, then the semi-Markov environment $Z(t)$ is called
\textit{Markov environment}.
\end{defn}

Another equivalent definition of Markov environment and
independent of Definition \ref{Defn-semiMI} is as follows.  Let
$Z(t)$ be a homogeneous Markov process with initial state
$\mathcal{E}(0)$ and transition probabilities $z_{l,m}\triangle
t+o(\triangle t)$ from the state $\mathcal{E}_l$ to the state
$\mathcal{E}_m$ ($l\neq m$) during a small time interval ($t$,
$t+\triangle t$), and there is probability $1-\sum_{m\neq
l}z_{l,m}\triangle t+o(\triangle t)$ to stay in the same state
$\mathcal{E}_l$ during the same time interval ($t$, $t+\triangle
t$).

In many papers on queueing theory, telecommunication systems,
inference of stochastic processes, statistics and other areas, the
above Markov environment is often associated with \textit{Markov
Modulated Poisson Process}. The phrases \textit{Markov
environment}, \textit{piece-wise deterministic environment} and
\textit{semi-Markov environment} are more appropriate in the
context of the present paper.

\subsection{Classification of bottleneck stations}
There was mentioned that $\mathcal{E}(t)$ denotes the state of the
semi-Markov environment in time $t$. For example, the equality
$\mathcal{E}(4)=\mathcal{E}_5$ means that in time $t=4$ the
associated Markov chain is in state $\mathcal{E}_5$. We also use
the following notation:
$\lambda_{j}(\mathcal{E}(t))=\lambda(\mathcal{E}(t))p_{j}(\mathcal{E}(t))$.

In the sequel we shall also use the notation $\lambda(t)$,
$p_{j}(t)$ and correspondingly, $\lambda_{j}(t)$. It is worth
noting, that the definition of $\lambda(t)$ and
$\lambda(\mathcal{E}(t))$ and consequently $\lambda_j(t)$ and
$\lambda_j(\mathcal{E}(t))$ have the different meaning. In general
$\lambda(\mathcal{E}(t))$ $\neq\lambda(t)$ and
$\lambda_j(\mathcal{E}(t))$ $\neq\lambda_j(t)$. $\lambda_{j}(t)$
and $\lambda_j(\mathcal{E}(t))$ both are random parameters, and
$\lambda_{j}(t)=\lambda_{j}(t,\omega)$ where $\omega\in\Omega$,
while $\lambda_j(\mathcal{E}(t))=\lambda_j(t,\omega_Z)$ where
$\omega_Z\in\Omega_Z$. (In the sequel this dependence upon
$\omega_Z$ will be always implied, but shown explicitly only in
the cases where it is necessary.)

More specific explanation of the above difference is as follows.
If at the initial time moment $t=0$ all of units are in the server
station, then the input rate to the $j$th client station is
$\lambda_j(0)N$. ($\lambda_j(0)$ is the individual rate of each
unit arriving to the client station $j$, and therefore the rate
between arrivals is $\lambda_j(0)N$.) The time parameter 0 in
parentheses is associated with the state of Markov environment in
time 0. Specifically, if at the initial time moment $t=0$ all of
units are in the server station, then
$\lambda_j(0)=\lambda_j(\mathcal{E}(0))$. If the network is
considered without Markov or semi-Markov environment, then that
initial arrival rate is $\lambda_jN$ which is associated with the
individual rate $\lambda_j$ of each unit of the server station. If
at the initial time moment there are $\alpha N$ units in the
server station, $\alpha<1$, then for arrival process to any client
station $j$ of a standard network (without random environment) we
also use the notation $\lambda_j(0)$. However, the meaning of
$\lambda_j(0)$ is not longer the individual rate of each unit at
time $t=0$ arriving to the station $j$. More specifically,
$\lambda_j(0)=\lambda_j\alpha$, where $\lambda_j$ is the
individual rate of each unit at time $t=0$ arriving to the station
$j$. The meaning of $\lambda_j(t)$ is similar. The only difference
that it is said about an arbitrary time $t$. For example, if there
are $\alpha(t)N$ units in the server station in time $t$,
$\alpha(t)<1$, then $\lambda_j(t)=\lambda_j\alpha(t)$. Resuming
the above, $\lambda_j(\mathcal{E}(t))$ is associated with
individual service rate of each unit, while $\lambda_j(t)$ is a
recalculated (relative) rate depending on the state of the
queue-length processes in time $t$, as explained above.

Regarding the departure (service) rates our assumption in the
paper is as follows. The departure rate of the $j$th client
station is assumed to be independent of semi-Markov environment as
well as of input rates $\lambda_j(t)$. Specifically,

\smallskip
 $\bullet$ it is assumed that the the expectation of
service (inter-departure) time in the $j$th client station is
$\mathbb{E}\xi_{j,N,l}=\frac{1}{\mu_jN}$ for all $l$=1,2,\ldots

\smallskip
Therefore, if at the initial time moment $t=0$ all of units are in
the server station, then the load parameter of the $j$th client
station is $\rho_j(0)=\frac{\lambda_j(0)}{\mu_j}$, and in the case
where there is no semi-Markov environment,
$\rho_j(0)=\rho_j=\frac{\lambda_j}{\mu_j}$.

In the case of network in Markov environment, the meaning of the
notation $\lambda_j(t)$ is the same as well. If there are
$N\alpha(t)$ units in the server station in time $t$,
$\alpha(t)<1$, and the rate of arrival of a unit from the server
to the client station $j$ is $\lambda_j(\mathcal{E}(t))$, then
$\lambda_j(t)$ = $\lambda_j(\mathcal{E}(t))\alpha(t)$. Then the
load of the $j$th client station in time $t$ is
$\rho_j(t)=\frac{\lambda_j(t)}{\mu_j}$.

Now introduce necessary definitions. The first two definitions are
related to both standard client/server networks and client/server
networks in semi-Markov environment.

\begin{defn}\label{defn1}
The client station $j$ is called \textit{locally non-bottleneck}
in time $t$ if $\rho_j(t)<1$. Otherwise, the $j$th client station
is called \textit{locally bottleneck} in time $t$. A client
station locally (non-) bottleneck in time 0 will be also called
\textit{initially (non-) bottleneck}.
\end{defn}

\begin{defn}\label{defn2}
The client station $j$ is called \textit{non-bottleneck} in time
interval [$t_1, t_2$] if it is locally non-bottleneck in all
points of this interval. Otherwise, if there is a point $t^*\in
[t_1, t_2$] such that $\rho_j(t^*)\geq1$, then the client station
is called \textit{bottleneck} in time interval [$t_1, t_2$]. A
client station is called \textit{(non-) bottleneck} if it is
(non-) bottleneck for all $t$.
\end{defn}

\smallskip
The special definition for standard client/server networks
(without semi-Markov environment) is as follows.

\begin{defn}\label{defn3}
A client station is called \textit{absolutely non-bottleneck} if
it is a locally non-bottleneck station at the moment when all of
units are in the server station. Otherwise, a client station is
called \textit{absolutely bottleneck}.
\end{defn}

Clearly, that absolutely non-bottleneck client station $j$ is a
non-bottleneck client station, because then for all $t$ we have
$\lambda_j(t)<\mu_j$. In the next section we prove that absolutely
bottleneck client station is a bottleneck client station as well.
That is, if a client station is currently locally bottleneck
station, then it never can become a locally non-bottleneck. This
means that the client station is forever bottleneck.

Definition \ref{defn3} can be extended to network stations in
semi-Markov environment for an arbitrary time $t$. Specifically,
we have the following definition.

\begin{defn}\label{defn4}
The client station $j$ of a network in random environment is
called absolutely non-bottleneck in time $t$ if
$\lambda_j(\mathcal{E}(t))<\mu_j$. In other words, the client
station $j$ is absolutely non-bottleneck in time $t$ if in that
time $t$ the network belongs to some state $\mathcal{E}_i$ of the
environment, i.e. $\mathcal{E}(t)$ = $\mathcal{E}_i$, and in this
state $\lambda_j(\mathcal{E}_i)<\mu_j$. Otherwise, this client
station is called absolutely bottleneck in time $t$.
\end{defn}

The last notion enables us to judge on the behavior of client
stations in random intervals [$\sigma_i, \sigma_{i+1}$), where the
network is in given state $\mathcal{E}_l$. Recall that $\sigma_i$
is a time instant when the state of semi-Markov environment is
changed. So, during the random interval [$\sigma_i, \sigma_{i+1}$)
the network is in a fixed state of the semi-Markov environment.

\smallskip

\section{Queue-length processes in the client stations of networks
with semi-Markov environment}\label{Semi-Martingale} Consider a
client station $j$ ($j=1,2,\ldots,k$). Let $Q_{j,N}(t)$ denote a
queue-length there in time $t$. Assume that at the initial time
instant $t=0$, all of units are in the server station, i.e.
$Q_{j,N}(0)=0$ for all $j=1,2,\ldots,k$. This is the simplest
case, and we start from its study. For $t>0$,
\begin{equation}\label{QL1}
Q_{j,N}(t)=A_{j,N}(t)-D_{j,N}(t),
\end{equation}
where $A_{j,N}(t)$ is the arrival process to client station $j$,
and $D_{j,N}(t)$ is the departure process from that client station
$j$. The equation for departure process is as follows. Let
\begin{equation*}\label{QL2}
S_{j,N}(t) = \sum_{i=1}^\infty
\mathbb{I}\left\{\sum_{l=1}^i\xi_{j,N,l}\leq t\right\}, \ \
j=1,2,\ldots,k
\end{equation*}
be a point process associated with consecutive departures from the
$j$th client station. Then,
\begin{equation}\label{QL3}
\begin{aligned}
D_{j,N}(t)&=\int_0^t\mathbb{I}\{Q_{j,N}(s-)>0\}\mbox{d}S_{j,N}(s)\\
&=S_{j,N}(t)-\int_0^t\mathbb{I}\{Q_{j,N}(s-)=0\}\mbox{d}S_{j,N}(s),\\
& \ \ \ j=1,2,\ldots,k.
\end{aligned}
\end{equation}
The definition of the departure process given by \eqref{QL3} is as
in \cite{Abramov 2000}. However, the construction of arrival
process is more difficult.

Specifically,
\begin{equation}\label{QL5}
A_{j,N}(t)=\int_0^t\sum_{i=1}^{N}\mathbb{I}\left\{N-\sum_{l=1}^k
Q_{l,N}(s-)\geq i\right\}\mbox{d}\pi_{j,i}(s,\omega_Z).
\end{equation}
The processes $\{\pi_{j,i}(s,\omega_Z)\}$, $i=1,2,\ldots,N$,
appearing in relation \eqref{QL5} are a collection of
conditionally independent Poisson processes with parameters
depending on $\omega_Z$. This means the following. Assume that for
a given realization $\omega_Z$ we have a sequence:
0$<\sigma_1<\sigma_2<$\ldots, and for some $n$, $\sigma_n\leq
s<\sigma_{n+1}$. Assume that
$\mathcal{E}(\sigma_n)=\mathcal{E}_l$. Then
$\{\pi_{j,i}(s,\omega_Z)\}$, $i=1,2,\ldots,N$, is the sequence of
Poisson processes with the parameter $\lambda_j(\mathcal{E}_l)$.
That is the rate of Poisson process depends on the state of
semi-Markov environment in time $s$. This is just the main
difference between the consideration of \cite{Abramov 2000}, where
$\lambda_j$, $j=1,2,\ldots,k$, were non-random constants.

Relation \eqref{QL5} can be then rewritten
\begin{equation}\label{QL6}
A_{j,N}(t)=\sum_{v=1}^\infty\int_{t\wedge\sigma_{v-1}}^{t\wedge\sigma_v}
\sum_{i=1}^{N}\mathbb{I}\left\{N-\sum_{l=1}^k Q_{l,N}(s-)\geq
i\right\}\mbox{d}\pi_{j,i,v}(s),
\end{equation}
where $\pi_{j,i,v}(s)$ is an associated sequence of
(conditionally) independent Poisson process with parameter
depending on the state of the semi-Markov environment
$\mathcal{E}(\sigma_{v-1})$. (Here in \eqref{QL6} and later we use
the standard notation for a minimum of two numbers: $a\wedge
b\equiv\min(a,b)$.)

Relations \eqref{QL1} and \eqref{QL3} enable us to write the
following representation for the queue-length process
$Q_{j,N}(t)$:
\begin{equation}\label{QL7}
Q_{j,N}(t)=A_{j,N}(t)-S_{j,N}(t)+\int_0^t\mathbb{I}\{Q_{j,N}(s-)=0\}\mbox{d}S_{j,N}(s).
\end{equation}
This implies that $Q_{j,N}(t)$ is the normal reflection of the
process
\begin{equation}\label{QL8}
X_{j,N}(t)=A_{j,N}(t)-S_{j,N}(t), \ \ X_{j,N}(0)=0
\end{equation}
at zero. More accurately, $Q_{j,N}(t)$ is a nonnegative solution
of the Skorokhod problem (see \cite{Skorokhod 1961} as well as
\cite{Anulova and Liptser 1990}, \cite{Ramanan 2006}, \cite{Tanaka
1979} and others) on the normal reflection of the process
$X_{j,N}(t)$ at zero (for the detailed arguments see \cite{Kogan
and Liptser 1993}). Recall that according to the Skorokhod
problem,
\begin{equation}\label{QL9}
\int_0^t\mathbb{I}\{Q_{j,N}(s-)=0\}\mbox{d}S_{j,N}(s)=-\inf_{s\leq
t}X_{j,N}(s),
\end{equation}
and $Q_{j,N}(t)$ has the representation
\begin{equation}\label{QL10}
Q_{j,N}(t)=X_{j,N}(t)-\inf_{s\leq t}X_{j,N}(s).
\end{equation}
In the sequel, it is convenient to use the notation:
$$\Phi_t(X)=X(t)-\inf_{s\leq t}X(s)$$ for any c\'adl\'ag function
$X(t)$ satisfying $X(0)=0$ (see e.g. \cite{Kogan and Liptser
1993}). According to this notation, \eqref{QL10} can be rewritten
$Q_{j,N}(t)=\Phi_t(X_{j,N})$.

Next, we take into account that the process $A_{j,N}(t)$ is a
semimartingale adapted with respect to the filtration
$\mathcal{F}_t$. Let $\widehat{A}_{j,N}(t)$ denote the compensator
of $A_{j,N}(t)$ and $M_{A_{j,N}}(t)$ denote the square integrable
martingale of $A_{j,N}(t)$ in the Doob-Meyer semimartingale
decomposition: $A_{j,N}(t)$ = $\widehat{A}_{j,N}(t)$ +
$M_{A_{j,N}}(t)$, $j=1,2,\ldots,k$. Then, the process $X_{j,N}(t)$
given by \eqref{QL8} can be represented
\begin{equation}\label{QL11}
X_{j,N}(t)=\widehat{A}_{j,N}(t)-S_{j,N}(t)+M_{A_{j,N}}(t),
\end{equation}
where
\begin{equation}\label{QL12}
\widehat{A}_{j,N}(t)=\sum_{l=1}^\infty\int_{\sigma_{l-1}\wedge
t}^{\sigma_{l}\wedge
t}\lambda_j(\mathcal{E}(\sigma_{l-1}))\left\{N-\sum_{i=1}^k
Q_{i,N}(s)\right\}\mbox{d}s.
\end{equation}
The details for last formula \eqref{QL12} can be obtained from
\cite{Dellacherie 1972} or \cite{Liptser and Shiryayev 1989},
Theorem 1.6.1.

Let us now study asymptotic properties of the normalized
queue-lengths $q_{j,N}(t)$ = $\frac{Q_{j,N}(t)}{N}$ as
$N\to\infty$. For normalized processes we will use small Latin
letters. For example, $x_{j,N}(t)=\frac{X_{j,N}(t)}{N}$,
$\widehat{a}_{j,N}(t)=\frac{\widehat{A}_{j,N}(t)}{N}$, and so on.

Then, from \eqref{QL11} we have
\begin{equation}\label{QL13}
x_{j,N}(t)=\widehat{a}_{j,N}(t)-s_{j,N}(t)+m_{A_{j,N}}(t), \
j=1,2,\ldots,k.
\end{equation}
Let us derive a relation for
$\mathbb{P}^{\_}\lim_{N\to\infty}x_{j,N}(t)$.
($\mathbb{P}^{\_}\lim$ denotes the limit in probability.)

From Lenglart-Rebolledo inequality (e.g. Liptser and Shiryayev
\cite{Liptser and Shiryayev 1989}), we have:
\begin{equation}\label{QL14}
\begin{aligned}
\mathbb{P}\left\{\sup_{0\leq s\leq
t}|m_{A_{j,N}}(t)|>\delta\right\}&=\mathbb{P}\left\{\sup_{0\leq
s\leq t}|A_{j,N}(s)-\widehat{A}_{j,N}(s)|>\delta N\right\}\\
&\leq\frac{\epsilon}{\delta^2}+\mathbb{P}\{\widehat{A}_{j,N}(t)>\epsilon
N^2\},
\end{aligned}
\end{equation}
where $\epsilon$ = $\epsilon(N)$ vanishes such that $\epsilon
N\to\infty$. Then, by virtue of \eqref{QL12}, the term
$\mathbb{P}\{\widehat{A}_{j,N}(t)>\epsilon N^2\}$ =
$\mathbb{P}\{\widehat{a}_{j,N}(t)>\epsilon N\}$ vanishes as well,
and for any small $\delta>0$ the fraction
$\frac{\epsilon}{\delta^2}$ vanishes. Therefore,
\begin{equation}\label{QL15}
\mathbb{P}^{\_}\lim_{N\to\infty}m_{A_{j,N}}(t) =0
\end{equation}
for all $j=1,2,\ldots,k$ and $t\geq0$.

Next, according to the assumption above,
$\mathbb{P}^{\_}\lim_{N\to\infty}s_{j,N}(t)=\mu_jt$.

Therefore,
\begin{equation}\label{QL16}
\begin{aligned}
x_j(t)&=\mathbb{P}^{\_}\lim_{N\to\infty}x_{j,N}(t)\\
&=\mathbb{P}^{\_}\lim_{N\to\infty}\widehat{a}_{j,N}(t)-\mu_jt\\
&=\sum_{l=1}^\infty\int_{\sigma_{l-1}\wedge t}^{\sigma_{l}\wedge
t}\left[\lambda_j(\mathcal{E}(\sigma_{l-1}))\left\{1-\sum_{i=1}^k
\Phi_s(x_i)\right\}-\mu_j\right]\mbox{d}s.
\end{aligned}
\end{equation}
Representation \eqref{QL16} is the extension of the similar result
of \cite{Abramov 2000} for queue-length processes in client
stations of standard client/server networks. Recall that the
representation obtained in \cite{Abramov 2000} and \cite{Kogan and
Liptser 1993} is
\begin{equation}\label{QL16+}
x_j(t)=\int_{0}^{ t}\left[\lambda_j\left\{1-\sum_{i=1}^k
\Phi_s(x_i)\right\}-\mu_j\right]\mbox{d}s.
\end{equation}

As we can see, the representations given by \eqref{QL16} and
\eqref{QL16+} are similar. The only difference is in the presence
of infinite sum containing the integrals with random upper and
lower bounds, and \eqref{QL16+} is a particular case of
\eqref{QL16}.

\section{Bottleneck analysis in the case of a piece-wise constant environment}
\label{Piece-wise deterministic}

In this section we discuss the behavior of the queue-length
processes in a large closed client/server network in the case of
piece-wise constant environment. For this case relation
\eqref{QL16} reduces to
\begin{eqnarray}\label{PCE}
x_j(t) &=&\sum_{l=1}^r\int_{\sigma_{l-1}}^{\sigma_{l}}
\left[\lambda_j\left\{1-\sum_{i=1}^k
\Phi_s(x_i)\right\}-\mu_j\right]\mbox{d}s\\
&&+\int_{\sigma_{r}}^{t} \left[\lambda_j\left\{1-\sum_{i=1}^k
\Phi_s(x_i)\right\}-\mu_j\right]\mbox{d}s,\nonumber
\end{eqnarray}
where $r$ is the number of the state changes before time $t$.

Therefore, the bottleneck analysis of the network in piece-wise
constant environment reduces to the analysis of the network in
traditional formulation (without random environment) in some given
time intervals such as [$\sigma_{l-1}, \sigma_l$),
$l$=1,2,\ldots,$r$, or [$\sigma_r, t$).

Notice, that the bottleneck analysis of the Markovian
client/server model has been originally studied by Kogan and
Liptser \cite{Kogan and Liptser 1993}. Then these results were
extended for the case of autonomous service mechanism in client
stations in \cite{Abramov 2000}. However, the results obtained in
both of these papers are related to a single special case and are
not enough for the purpose of our performance analysis. Therefore
we will study all possible cases including the behaviour of the
network under different initial lengths of queues in client
stations and several absolutely bottleneck and absolutely
non-bottleneck client stations.

In \cite{Kogan and Liptser 1993} and \cite{Abramov 2000} there has
only been considered the case where in the initial time moment
$t=0$ all of units are in the server station (i.e. the client
stations all are empty) and only one (the $k$th) client station is
a bottleneck station. (In this particular case the notions of
(non-)bottleneck and absolutely (non-)bottleneck client station
coincide.) Specifically, there has been proved the following
result in \cite{Abramov 2000}.

\begin{prop}\label{prop1}
Let $S_{j,N}^*(t)=\inf\{s>0: S_{j,N}(s)=S_{j,N}(t)\}$, and
$Q_{j,N}(t)$ denotes the queue-length in the $j$th client station
in time $t$. Under the assumption that the $k$th client section is
bottleneck, for $j=1,2,\ldots,k-1$ and for any $t>0$ we have:
\begin{eqnarray}
&&\lim_{N\to\infty}\mathbb{P}\{Q_{j,N}[S_{j,N}^*(t)]=0\}=1-\rho_j(t),\label{BA1-1}\\
&&\lim_{N\to\infty}\int_0^t\rho_j(s)\mathbb{P}\{Q_{j,N}(s)=l\}\mbox{d}s\nonumber\\
&&=\lim_{N\to\infty}\int_0^t
\mathbb{P}\{Q_{j,N}[S_{j,N}^*(s)]=l+1\}\mbox{d}s,\label{BA1-2}\\
&&\ \ \ l=0,1,\ldots,\nonumber
\end{eqnarray}
where
\begin{eqnarray}
\rho_j(t)&=&\rho_j(0)[1-q(t)],\label{BA2-1}\\
q(t)&=&\left(1-\frac{1}{\rho_k(0)}\right)(1-\mbox{e}^{-\rho_k(0)\mu_kt}).\label{BA2-2}
\end{eqnarray}
\end{prop}

In the case of Markovian network, relations \eqref{BA1-1},
\eqref{BA1-2} can be resolved explicitly, and queue-length
distribution in the non-bottleneck stations is the time-dependent
geometric distribution \cite{Abramov 2000},  \cite{Kogan and
Liptser 1993}.

The meaning of $1-q(t)$ in \eqref{BA2-1} is the limiting fraction
of units remaining at the server station in time $t$ as
$N\to\infty$. For example, in the case where $\rho_k(0)=1$ this
fraction remains the same at any time $t$ as initially, that is,
as $N$ large, the number of units in the server station remains
asymptotically equivalent to $N$. However, if $\rho_k(0)>1$, then
the number of units in the server station in time $t$ is
asymptotically equivalent to $N[1-q(t)]$. Then, the number of
units remaining in bottleneck station in time $t$ is
asymptotically equivalent to $Nq(t)$.

Let us study various cases of the client/server network with
bottlenecks. These cases will be studied in order of increasing
complexity. We possibly shorten  the proofs.

Assuming that the initial condition of the network is the same as
in Proposition \ref{prop1} (i.e. at the initial time $t=0$ all of
units are at the server station), let us study the case where the
client stations 1,2,\ldots,$k_0$ are non-bottleneck, while the
rest client stations $k_0+1,\ldots,k$ are bottleneck. In this case
we have the following result.

\begin{prop}\label{prop2} Under the
assumption that the client sections $k_0+1,\ldots,k$ are
bottleneck, for $j=1,2,\ldots,k_0$ and for any $t>0$ we have:
\begin{eqnarray}
&&\lim_{N\to\infty}\mathbb{P}\{Q_{j,N}[S_{j,N}^*(t)]=0\}=1-\rho_j(t),\label{BA3-1}\\
&&\lim_{N\to\infty}\int_0^t\rho_j(s)\mathbb{P}\{Q_{j,N}(s)=l\}\mbox{d}s\nonumber\\
&&=\lim_{N\to\infty}\int_0^t
\mathbb{P}\{Q_{j,N}[S_{j,N}^*(s)]=l+1\}\mbox{d}s,\label{BA3-2}\\
&&\ \ \ l=0,1,\ldots,\nonumber
\end{eqnarray}
where
\begin{eqnarray}
\rho_j(t)&=&\rho_j(0)[1-q(t)],\label{BA4-1}\\
q(t)&=&\left(1-\frac{\sum_{v=k_0+1}^k\mu_v}{\sum_{v=k_0+1}^k\lambda_v}\right)
\left(1-\exp\left[-t\sum_{v=k_0+1}^k\lambda_v\right]\right).\label{BA4-2}
\end{eqnarray}
\end{prop}

Note, that equations \eqref{BA3-1} and \eqref{BA3-2} are the same
as \eqref{BA1-1} and \eqref{BA1-2}, and the main difference
between Proposition \ref{prop1} and Proposition \ref{prop2} is
only in the expression for $q(t)$. The difference between
expressions \eqref{BA2-1}, \eqref{BA2-2} and \eqref{BA4-1},
\eqref{BA4-2} can be easily explained in the framework of the
proof of Proposition \ref{prop2}.

Considering all of the bottleneck client stations as a separate
subsystem, one can notice that the arrival rate to this subsystem
is $\sum_{v=k_0+1}^k\lambda_v$, and the service rate (the sum of
reciprocals of the expected inter-departure times) is
$\sum_{v=k_0+1}^k\mu_v$. This subsystem can be thought as a
bottleneck station with the load
\begin{equation*}
\label{BA5}\frac{\sum_{v=k_0+1}^k\lambda_v}{\sum_{v=k_0+1}^k\mu_v}.
\end{equation*}
To be specific, note the method of Sections 2, 3, and 4 of
\cite{Abramov 2000} leads to the same equations, and the equations
for normalized queue-lengths in the bottleneck stations all are an
elementary extension of the case considered in \cite{Abramov
2000}.

Therefore Proposition \ref{prop2} is an elementary extension of
Proposition \ref{prop1} and its proof is the same as in
\cite{Abramov 2000} or \cite{Kogan and Liptser 1993}. Notice, that
the number of units remaining at the bottleneck station $v$,
$v=k_0+1,\ldots,k$ in time $t$ is asymptotically equal to
\begin{equation}
\label{BA6}
Nq_v(t)=N\left([\lambda_v-\mu_v]t-\lambda_v\int_0^tq(s)\mbox{d}s\right).
\end{equation}

Recall the main elements of the known proof for the representation
$q(t)$ given by \eqref{BA2-2}, and consequently explain the proof
of \eqref{BA6}. We use the notation similar to that of the earlier
papers \cite{Abramov 2000} and \cite{Kogan and Liptser 1993}.

The difference between arrival and departure processes in the
$j$th client station is denoted $X_{j,N}(t)=A_{j,N}(t)-S_{j,N}(t)$
and its normalization $x_{j,N}(t)=\frac{X_{j,N}(t)}{N}$. Let
$x_j(t)$, $j=1,2,\ldots,k$, denote the limit in probability of
$x_{j,N}(t)$, as $N\to\infty$. The queue-length in the $j$th
client station in time $t$ is denoted $Q_{j,N}(t)$, its
normalization is denoted $q_{j,N}(t)=\frac{Q_{j,N}(t)}{N}$, and
the limit of $q_{j,N}(t)$ in probability, as $N\to\infty$, is
denoted $q_j(t)$.

Next, let $\Phi_t(X)=X(t)-\inf_{s\leq t}X(s)$ for any c\'adl\'ag
function $X$ satisfying the condition $X(0)=0$. The functional
$\Phi_t(X)$ has been introduced in \cite{Kogan and Liptser 1993}.
It characterizes a solution of the Skorokhod problem on normal
reflection at zero. Therefore $\Phi_t(x_j)$, being the functional
$\Phi$ applied to the function $x_j(t)$, describes the dynamic of
normalized queue-length in the $j$th client station under the
``usual" initial conditions given in Propositions \ref{prop1} and
\ref{prop2}. Under these ``usual" conditions, the functions
$x_j(t)$, $j=1,2,\ldots,k$ satisfy the system of equations:
\begin{equation}\label{BA7}
x_j(t)=\int_0^t\left\{\lambda_j\left[1-\sum_{l=1}^k\Phi_s(x_l)\right]-\mu_j\right\}\mbox{d}s,
\ \ \ j=1,2,\ldots,k.
\end{equation}
Note, that the normalized functions $x_j(t)$, $j=1,2,\ldots,k$ are
usual (non-random) continuous functions, and \eqref{BA7}
characterizes a usual system of linear differential equations.

The statement of Proposition \ref{prop1} is based on the solution
of the system of these equations. (It is proved in \cite{Kogan and
Liptser 1993} that there is a unique solution of the system of
equations \eqref{BA7}.) More specifically, in the case where the
node $k$ is bottleneck, $\inf_{s\leq t} x_k(s)=x_k(0)=0$, and we
therefore have $\Phi_s(x_k)=x_k(s)$. In other words, $\Phi_s(x_k)$
can be replaced by $x_k(s)$ in these equations of \eqref{BA7}.

The solution of the system of equations \eqref{BA7} is $x_j(t)$=0
for the non-bottleneck stations $j=1,2,\ldots,k-1$, and
$x_k(t)=q(t)$ for the bottleneck station $j=k$, where $q(t)$ is
given by \eqref{BA2-2}. In the case of several bottleneck
stations, we write the similar equation for the dynamic of the
normalized cumulated queue-length process in all of the bottleneck
stations. Specifically,
\begin{equation*}
\sum_{v=k_0+1}^k
x_v(t)=\int_0^t\left\{\sum_{v=k_0+1}^k\lambda_v\left[1-\sum_{l=k_0+1}^k
x_l(s)\right]-\sum_{v=k_0+1}^k\mu_v\right\}\mbox{d}s,
\end{equation*}
gives solution \eqref{BA4-2} for $q(t)$ in the similar statement
of Proposition \ref{prop2}. Then, the solution of system
\eqref{BA7} for the bottleneck client stations
$v$=$k_0+1$,$k_0+2$\ldots, $k$ is given by
\begin{equation*}
x_v(t)=q_v(t)=(\lambda_v-\mu_v)t-\lambda_v\int_0^t q(s)\mbox{d}s,
\end{equation*}
where $q(t)$ is defined by \eqref{BA4-2}, and the queue-lengths in
the bottleneck stations are asymptotically evaluated by relation
\eqref{BA6}.

Now, we discuss the behavior of the network, in which the client
stations are not initially empty. This is the next step of the
extension of the original result of Proposition \ref{prop1}. The
analysis of cases related to initially not empty queues is much
more difficult. Therefore, we start from the simplest case of the
network containing only one client station, i.e. $k=1$.

Let $\beta_1\leq1$ be a positive real number, and let us assume
that the initial number of units in this client station is
asymptotically equivalent to $N\beta_1$. Consider the following
two cases: (i) the client station is initially bottleneck, i.e.
$\lambda_1(0)\geq\mu_1$, and (ii) the client station is initially
non-bottleneck, i.e. $\lambda_1(0)<\mu_1$.

Case (i) is relatively simple. It is a simple extension of the
cases considered above. Specifically, we have the following system
of equations:
\begin{equation}\label{BA8}
\begin{aligned}
&x_1(t)=\beta_1+(1-\beta_1)z_1(t),\\
&z_1(t)=\int_0^t\left\{\lambda_1(0)\left[1-
z_1(s)\right]-\mu_1\right\}\mbox{d}s.
\end{aligned}
\end{equation}
From \eqref{BA8} we have the following solution:
\begin{equation}\label{BA9}
x_1(t)=\beta_1+(1-\beta_1)\left(\frac{\lambda_1(0)-\mu_1}{\lambda_1(0)}\right)\left(1-\mbox{e}^{-\lambda_1(0)t}\right).
\end{equation}
The normalized queue-length $q_1(t)$ in this client station is
$q_1(t)=x_1(t)$.

Case (i) can be easily extended to a more general case of $k$
initially bottleneck client stations. Let $\beta_j$,
$j=1,2,\ldots,k$, denote nonnegative real numbers, and $\beta_1$ +
$\beta_2$ + \ldots + $\beta_k\leq1$. Assume then that the initial
number of units in the $j$th client station is asymptotically
equivalent to $N\beta_j$. Let $q(t)$ denote the cumulated
normalized queue-length in all of client stations.

\begin{prop}\label{prop3}
Assume that all client stations are initially bottleneck, and the
initial queue-lengths in client stations are asymptotically
equivalent to $N\beta_1$, $N\beta_2$,\ldots, $N\beta_k$
correspondingly ($\beta_1 + \beta_2 +\ldots+ \beta_k\leq1$), as
$N\to\infty$. Then,
\begin{equation}\label{BA10}
\begin{aligned}
&q(t)=\sum_{j=1}^k\beta_j+\left(1-\sum_{j=1}^k\beta_j\right)r(t),\\
&r(t)=\left(\frac{\sum_{j=1}^k(\lambda_j(0)-\mu_j)}{\sum_{j=1}^k\lambda_j(0)}\right)
\left(1-\exp\left[ -t\sum_{j=1}^k\lambda_j(0)\right]\right),
\end{aligned}
\end{equation}
and the normalized queue-length in the $j$th client station is
defined as
\begin{equation}\label{BA11}
q_j(t)=\beta_j+\left(1-\sum_{j=1}^k\beta_j\right)\left([\lambda_j(0)-\mu_j]t-\lambda_j(0)\int_0^t
r(s)\mbox{d}s\right).
\end{equation}
\end{prop}

Let us now discuss case (ii). This case is also described by
system of equations \eqref{BA8}, and the dynamic of the process
$x_1(t)$ is therefore similar to the case considered above.
However, this case is more delicate. The client station is
initially non-bottleneck, i.e. $\lambda_1(0)<\mu_1$, and the
function $x_1(t)$ is therefore decreasing in the right side of 0.
According to the convention, the initial value of queue is
asymptotically equivalent to $\beta_1N$, and therefore
$\lambda_1(0)=(1-\beta_1)\lambda_1^*$. Then the meaning of
$\lambda_1^*$ is a maximally possible rate of units arriving from
the server station to client station, when all of units are in the
server station and the client station is empty. Then the client
station is absolutely non-bottleneck if $\lambda_1^*<\mu_1$, and
it is absolutely bottleneck if $\lambda_1^*\geq\mu_1$.

Consider first the case of an absolutely bottleneck station, i.e.
$\lambda_1^*\geq\mu_1$. In this case, from the solution given by
\eqref{BA9} we have
\begin{equation}\label{BA12}
\begin{aligned}
\lim_{t\to\infty}q_1(t)&=\beta_1+(1-\beta_1)\frac{\lambda_1(0)-\mu_1}{\lambda_1(0)}\\
&=\frac{\lambda_1(0)-\mu_1(1-\beta_1)}{\lambda_1(0)}\\
&=\frac{\lambda_1^*-\mu_1}{\lambda_1^*}.
\end{aligned}
\end{equation}
The meaning of the last result is the following. Let
$N\overline{\beta_1}$ be an asymptotic value of the queue-length
in the client station, when this queue station (at the first time)
becomes a locally bottleneck station. Then, according to
\eqref{BA12} the normalized queue-length in the client station
approaches to this level
$\overline{\beta_1}=\frac{\lambda_1^*-\mu_1}{\lambda_1^*}$ as
$t\to\infty$. Notice, that the same level for normalized
queue-length is achieved for an initially bottleneck station in
case (i). The result of \eqref{BA12} is also supported by results
\eqref{BA2-1} and \eqref{BA2-2} of Proposition \ref{prop1}. Thus,
for a bottleneck station, the same level is asymptotically
achieved independently of an initial queue-length. For this reason
for any absolutely bottleneck station in which
$\lambda_1^*\geq\mu_1$ we do not distinguish between two cases
$\lambda_j(0)<\mu_j$ and $\lambda_j(0)\geq\mu_j$, and absolutely
bottleneck client station is always a bottleneck station.

Consider now the case of an absolutely non-bottleneck client
station (and therefore non-bottleneck client station) where
$\lambda_1^*<\mu_1$. Then, according to the same calculation as in
\eqref{BA12} we have
\begin{equation}\label{BA13}
\lim_{t\to\infty}x_1(t)=\frac{\lambda_1^*-\mu_1}{\lambda_1^*}<0.
\end{equation}
Therefore, there exists the time instant $\tau_1$ when the
normalized queue-length becomes at the first time empty. For this
time instant we have the equation
\begin{equation}\label{BA14}
\tau_1=-\frac{1}{\lambda_1(0)}\log\left(1+\frac{\beta_1\lambda_1(0)}
{(1-\beta_1)(\lambda_1(0)-\mu_1)}\right).
\end{equation}

Let us extend the result of case (ii) for a network with $k$
initially non-bottleneck client stations, all satisfying the
condition $\lambda_j(0)<\mu_j$, $j=1,2,\ldots,k$. It is assumed
that the initial number of units in the $j$th client station is
asymptotically equivalent to $N\beta_j$ ($\beta_1$ + $\beta_2$
+\ldots+ $\beta_k\leq1$) as $N\to\infty$. Assume also that the
first $k_0$ client stations are (absolutely) non-bottleneck, i.e.
$\lambda_j^*$ = $\lambda_j$(0)(1$-\beta_1 - \beta_2 -\ldots-
\beta_k$) $<\mu_j$, $j=1,2,\ldots,k_0$, while the rest $k-k_0$
client stations are (absolutely) bottleneck, i.e. $\lambda_v^*$ =
$\lambda_v$(0)(1$-\beta_1 - \beta_2 -\ldots- \beta_k$)
$\geq\mu_v$, $v=k_0+1,k_0+2,\ldots,k$.

Similarly to \eqref{BA8} we have the following system of equations
\begin{equation}\label{BA15}
\begin{aligned}
&x(t)=\sum_{j=1}^k\beta_j+\left(1-\sum_{j=1}^k\beta_j\right)z(t),\\
&z(t)=\int_0^t\left\{(1-z(s))\sum_{j=1}^k\lambda_j(0)-\sum_{j=1}^k\mu_j\right\}\mbox{d}s,
\end{aligned}
\end{equation}
and for $z(t)$ we have the solution
\begin{equation}\label{BA16}
z(t)=\left(\frac{\sum_{j=1}^k(\lambda_j(0)-\mu_j)}{\sum_{j=1}^k\lambda_j(0)}\right)
\left(1-\exp\left[-t\sum_{j=1}^k\lambda_j(0)\right]\right).
\end{equation}

Then for $x_j(t)$, $j=1,2,\ldots,k$, we have the following
solutions:
\begin{equation}\label{BA17}
\begin{aligned}
x_j(t)&=\beta_j+\left(1-\sum_{j=1}^k\beta_j\right)\left((\lambda_j(0)-\mu_j)t-\lambda_j(0)\int_0^tz(s)\mbox{d}s\right).\\
\end{aligned}
\end{equation}
(Recall that
$x_j(t)=\mathbb{P}^{\_}\lim_{N\to\infty}\frac{A_{j,N}(t)-S_{j,N}(t)}{N}$,
$j=1,2,\ldots,k$.)

However, since the first $k_0$ client stations are non-bottleneck,
then the equality $q_j(t)=x_j(t)$ for the normalized queue-lengths
in client stations is valid only for the values $t$ of the
interval $0\leq t\leq\tau_1$, where the value $\tau_1$ can be
found from \eqref{BA17} as
\begin{equation}\label{BA18}
\tau_1=\min_{1\leq j\leq k_0}\inf\{t: x_j(t)\leq0\}.
\end{equation}
Let $j_0$ = $\arg\min_{1\leq j\leq k_0}\inf\{t: x_j(t)\leq0\}$.
Then the normalized queue-length process $q_{j_0}(t)$ is as
follows. For $0\leq t\leq\tau_1$, $q_{j_0}(t)=x_{j_0}(t)$, and for
$t\geq\tau_1$ it satisfies the equation
\begin{equation}\label{BA19}
x_{j_0}(t)=\int_{0}^{t-\tau_1}\left[\lambda(\tau_1)\left(1-\sum_{j\neq
j_0}x_j(s)-\Phi_s(x_{j_0})\right)-\mu_{j_0}\right]\mbox{d}s,
\end{equation}
where $\Phi_s(x_{j_0})=x_{j_0}(s)-\inf_{0\leq u\leq s}x_{j_0}(u)$.

Together with \eqref{BA19} for all remaining $j$ =1,2,\ldots,
$j_0-1$, $j_0+1$,\ldots, $k$ and $t\geq\tau_1$ we have:
\begin{equation}\label{BA20}
\begin{aligned}
&x(t)=\sum_{j\neq j_0} x_j(\tau_1)+\left(1-\sum_{j\neq j_0} x_j(\tau_1)\right)z(t),\\
&z(t)=\int_{0}^{t-\tau_1}\left\{(1-z(s))\sum_{j\neq
j_0}\lambda_j(\tau_1)-\sum_{j\neq j_0}\mu_j\right\}\mbox{d}s,
\end{aligned}
\end{equation}
and similarly to \eqref{BA16} for $z(t)$ we have the solution
\begin{equation}\label{BA21}
z(t)=\left(\frac{\sum_{j\neq
j_0}(\lambda_j(\tau_1)-\mu_j)}{\sum_{j\neq j_0}\lambda_j(\tau_1)}
\right) \left(1-\exp\left[-(t-\tau_1)\sum_{j\neq
j_0}\lambda_j(\tau_1)\right]\right).
\end{equation}
Therefore for $j$ =1,2,\ldots, $j_0$-1, $j_0$+1,\ldots, $k$ and
$t\geq\tau_1$ we obtain
\begin{equation}\label{BA22}
\begin{aligned}
x_j(t)=&x_j(\tau_1)+\left(1-\sum_{j\neq
j_0}x_j(\tau_1)\right)\times\\
&\times\left((\lambda_j(\tau_1)-\mu_j)(t-\tau_1)-\lambda_j(\tau_1)
\int_0^{t-\tau_1}z(s)\mbox{d}s\right),
\end{aligned}
\end{equation}
and since $x_{j_0}(t)$ is nonnegative, for $j=j_0$ and
$t\geq\tau_1$ we obtain
\begin{equation}\label{BA23}
x_{j_0}(t)=0.
\end{equation}

Thus the dimension of the system is decreased by 1, and the
procedure can be repeated similarly. Specifically, again since the
client stations $j$ =1,2,\ldots, $j_0-1$, $j_0+1$,\ldots, $k_0$
all are non-bottleneck stations, then the equality $q_j(t)=x_j(t)$
for the normalized queue-length processes is valid for all
$t\leq\tau_2$, where $\tau_2$ is defined as
\begin{equation}\label{BA24}
\tau_2=\min_{1\leq j\neq j_0\leq k_0}\inf\{t: x_j(t)\leq0\}.
\end{equation}
Setting now $j_1=\arg\min_{1\leq j\neq j_0\leq k_0}\inf\{t:
x_j(t)\leq0\}$ one can continue this procedure to find
$\tau_3$,\ldots, $\tau_{k_0}$ and then to know the behaviour of
the queue-length processes in all of (non-bottleneck and
absolutely bottleneck) client stations for all $t$.

The considered extension of case (ii) is in fact the general case
in which amongst $k$ client stations there are $k_0$
non-bottleneck client stations, and the rest $k-k_0$ stations are
absolutely bottleneck (i.e. part of them can be initially
non-bottleneck), and all $k$ client stations are with arbitrarily
large initial queue-lengths.

\smallskip
The result can be formulated as follows.
\begin{thm}
\label{thm1} Assume that there are $k$ client stations, where
$k_0$ client stations (non necessarily the first ones) are
absolutely non-bottleneck, and the initial normalized queue-length
in all of these $k$ client stations are $\beta_j$ in limit as
$N\to\infty$ correspondingly ($j=1,2,\ldots,k$). Then there are
time instants $\tau_1\leq\tau_2\leq\ldots\leq\tau_{k_0}$ when the
normalized queue-lengths in these $k_0$ absolutely non-bottleneck
stations correspondingly achieves zero and remains then to stay at
zero. These points as well as the normalized queue-lengths in time
$t$ in the rest $k-k_0$ stations are defined according to Scheme
\ref{alg1} below.

\begin{alg}\label{alg1}
 Consider the system
\begin{equation}\label{BA25}
\begin{aligned}
x_j(t)&=x_j(0)+\left(1-\sum_{j=1}^k x_j(0)\right)
\left((\lambda_j(0)-\mu_j)t-\lambda_j(0)\int_0^tz(s)\mbox{d}s\right),\\
&j=1,2,\ldots,k,
\end{aligned}
\end{equation}
where $x_j(0)=\beta_j$, and
\begin{equation}\label{BA26}
z(t)=\left(\frac{\sum_{j=1}^k(\lambda_j(0)-\mu_j)}{\sum_{j=1}^k\lambda_j(0)}\right)
\left(1-\exp\left[-t\sum_{j=1}^k\lambda_j(0)\right]\right).
\end{equation}
Then,
\begin{equation*}
\tau_1=\min_{1\leq j\leq k}\inf\{t: x_j(t)\leq0\},
\end{equation*}
and the argument $t$ in \eqref{BA25} and \eqref{BA26} belongs to
the interval [0,$\tau_1$]. In this case the normalized
queue-length at the $j$th client station, $q_j(t)=x_j(t)$ for all
$j$=1,2,\ldots,$k$. Let $j_0$ = $\arg\min_{1\leq j\leq k}\inf\{t:
x_j(t)\leq0\}$. Then $q_{j_0}(t)$=0 for all $t\geq\tau_1$.

In the next step, we consider the system of $k-1$ equations, and
$t\geq\tau_1$ (the equation for $x_{j_0}(t)$, $t\geq\tau_1$, is
excluded):
\begin{equation}\label{BA27}
\begin{aligned}
x_j(t)=&x_j(\tau_1)+\left(1-\sum_{j\neq
j_0}x_j(\tau_1)\right)\times\\
&\times\left((\lambda_j(\tau_1)-\mu_j)(t-\tau_1)-\lambda_j(\tau_1)
\int_0^{t-\tau_1}z(s)\mbox{d}s\right),
\end{aligned}
\end{equation}
where
\begin{equation}\label{BA29}
z(t)=\left(\frac{\sum_{j\neq
j_0}(\lambda_j(\tau_1)-\mu_j)}{\sum_{j\neq j_0}\lambda_j(\tau_1)}
\right) \left(1-\exp\left[-t\sum_{j\neq
j_0}\lambda_j(\tau_1)\right]\right).
\end{equation}
Then,
\begin{equation*}
\begin{aligned}
 &\tau_2=\min_{\substack{1\leq j\leq k\\ j\neq j_0}}\inf\{t\geq\tau_1:
x_j(t)\leq0\},\\
&j_1=\arg\min_{\substack{1\leq j\leq k\\ j\neq
j_0}}\inf\{t\geq\tau_1: x_j(t)\leq0\},
\end{aligned}
\end{equation*}
and the argument $t$ in \eqref{BA27} and \eqref{BA29} belongs to
the interval [$\tau_1$,$\tau_2$], and subscript $j$ in these
equations is $j=1,2,\ldots,j_0-1,j_0+1,\ldots,k$. In this case
$q_j(t)=x_j(t)$, and $q_{j_1}(t)=0$ for all $t\geq\tau_2$.

Consequently excluding equations one-by-one, then the $l$th
instant $\tau_l$, $2\leq l\leq k_0$, is defined from the system
\begin{equation}\label{BA30}
\begin{aligned}
x_j(t)=&x_j(\tau_{l-1})+\left(1-\sum_{j\neq
j_0,j_1,\ldots,j_{l-2}}x_j(\tau_{l-1})\right)\times\\
&\times\left((\lambda_j(\tau_{l-1})-\mu_j)(t-\tau_{l-1})-\lambda_j(\tau_{l-1})
\int_0^{t-\tau_{l-1}}z(s)\mbox{d}s\right),
\end{aligned}
\end{equation}
where
\begin{equation}\label{BA32}
\begin{aligned}
z(t)=&\left(\frac{\sum_{j\neq
j_0,j_1,\ldots,j_{l-2}}(\lambda_j(\tau_{l-1})-\mu_j)}{\sum_{j\neq
j_0,j_1,\ldots,j_{l-2}}\lambda_j(\tau_{l-1})} \right)\times\\
&\times\left(1-\exp\left[-t\sum_{j\neq
j_0,j_1,\ldots,j_{l-2}}\lambda_j(\tau_{l-1})\right]\right).
\end{aligned}
\end{equation}
Then,
\begin{equation*}
\begin{aligned}
 &\tau_{l}=\min_{\substack{1\leq j\leq k\\ j\neq j_0,j_1,\ldots,j_{l-2}}}\inf\{t\geq\tau_{l-1}:
x_j(t)\leq0\},\\
&j_{l-1}=\arg\min_{\substack{1\leq j\leq k\\ j\neq
j_0,j_1,\ldots,j_{l-2}}}\inf\{t\geq\tau_{l-1}: x_j(t)\leq0\},
\end{aligned}
\end{equation*}
and the argument $t$ in \eqref{BA30} and \eqref{BA32} belongs to
the interval [$\tau_{l-1}$,$\tau_l$], and subscript $j$ in these
equations takes the values from 1 to $k$ but $j_0$, $j_1$,\ldots,
$j_{l-2}$. In this case the normalized queue-lengths
$q_j(t)=x_j(t)$, and $q_{j_{l-1}}(t)$=0 for all $t\geq\tau_{l}$.
\end{alg}
\end{thm}

Theorem \ref{thm1} containing Scheme \ref{alg1} is easily applied
to networks with piece-wise deterministic environment.
Specifically, the theorem is applied to each of the intervals
[$\sigma_{l-1}, \sigma_l$), $l=1,2,\ldots,r$ and [$\sigma_r, t$).

\section{The main theorem on queue-length processes in the client stations of networks
with semi-Markov environment}\label{Semi-Markov}

The statement of Theorem \ref{thm1} is easily adapted to
client/server networks in semi-Markov environment. The difference
between the approach to the piece-wise deterministic environment
and semi-Markov environment is only that the points $\sigma_l$ are
random, and general relation \eqref{QL16} rather than \eqref{PCE}
must be used.

For example, considering the first term of \eqref{QL16}, we have
\begin{equation}\label{QL16-1}
\begin{aligned}
&x_j(\sigma_1\wedge t)\\
&=\sum_{l=1}^\infty\int_{\sigma_{l-1}\wedge (\sigma_1\wedge
t)}^{\sigma_{l}\wedge (\sigma_1\wedge
t)}\left[\lambda_j(\mathcal{E}(\sigma_{l-1}))\left\{1-\sum_{i=1}^k
\Phi_s(x_j)\right\}-\mu_j\right]\mbox{d}s\\
&=\int_{0}^{\sigma_{1}\wedge
t}\left[\lambda_j(\mathcal{E}(0))\left\{1-\sum_{i=1}^k
\Phi_s(x_j)\right\}-\mu_j\right]\mbox{d}s,
\end{aligned}
\end{equation}
and the extension the above theory of Section \ref{Piece-wise
deterministic} is elementary. For example one can make the
assumption about arbitrary initial conditions in client stations
and arbitrary number of bottleneck stations. The same theory of
Section \ref{Piece-wise deterministic} can be developed for any
interval with random lower and upper bounds as well, such that
\begin{equation*}
\int_{\sigma_{l-1}\wedge t}^{\sigma_{l}\wedge
t}\left[\lambda_j(\mathcal{E}(\sigma_{l-1}))\left\{1-\sum_{i=1}^k
\Phi_s(x_j)\right\}-\mu_j\right]\mbox{d}s,
\end{equation*}
resulting in adaptation of the theory of Section \ref{Piece-wise
deterministic} to client/server networks in semi-Markov
environment.

The theorem below is the adaptation of Theorem \ref{thm1} to the
case of a random interval [$0$, $\sigma_{1}\wedge t$). (The
extension of Theorem \ref{thm1} to the other intervals
[$\sigma_{l-1}\wedge t$, $\sigma_{l}\wedge t$), $l$=1,2,\ldots, is
similar.)

\begin{thm}
\label{thm2} Assume that amongst $k$ client stations, there are
$k_0$ absolutely non-bottleneck in time 0, and the rest $k-k_0$
client stations are absolutely bottleneck in time 0. Assume that
the initial normalized queue-lengths in all of these $k$ client
stations are $\beta_j$ in limit as $N\to\infty$ correspondingly
($j=1,2,\ldots,k$). Then there are time instants
$\tau_1\leq\tau_2\leq\ldots\leq\tau_{k_0}$ which are defined
recurrently by Scheme \ref{alg1} of Theorem \ref{thm1}. We only
take into account the values $\tau_i$ satisfying the inequality
$\tau_i\leq\sigma_1\wedge t$, i.e. we set $\ell$ =
$\ell(\sigma_1)=\max\{i: \tau_i\leq\sigma_1\wedge t\}$. Then, the
only instants $\tau_1\leq$ $\tau_2\leq$\ldots $\leq\tau_\ell$ are
taken into account in this theorem. The main relation \eqref{BA25}
of Scheme \ref{alg1} now looks
\begin{equation}\label{QL17}
\begin{aligned}
x_j(\sigma_1\wedge t)=&x_j(0)+\left(1-\sum_{j=1}^k
x_j(0)\right)\times\\
&\times\left((\lambda_j(0)-\mu_j)(\sigma_1\wedge t)-\lambda_j(0)\int_0^{\sigma_1\wedge t}z(s)\mbox{d}s\right),\\
&j=1,2,\ldots,k,
\end{aligned}
\end{equation}
where $x_j(0)=\beta_j$, and
\begin{equation}\label{QL18}
z(t)=\left(\frac{\sum_{j=1}^k(\lambda_j(0)-\mu_j)}{\sum_{j=1}^k\lambda_j(0)}\right)
\left(1-\exp\left[-t\sum_{j=1}^k\lambda_j(0)\right]\right).
\end{equation}
The other relations of this theorem are defined similarly to the
corresponding relations of Scheme \ref{alg1} where only argument
$t$ is replaced by $\sigma_1\wedge t$ in the corresponding places.
Specifically, \eqref{BA27} now looks
\begin{equation*}\label{QL19}
\begin{aligned}
&x_j(\sigma_1\wedge t)=x_j(\tau_1)+\left(1-\sum_{j\neq
j_0}x_j(\tau_1)\right)\times\\
&\times\left((\lambda_j(\tau_1)-\mu_j)((\sigma_1\wedge
t)-\tau_1)-\lambda_j(\tau_1) \int_0^{(\sigma_1\wedge
t)-\tau_1}z(s)\mbox{d}s\right),
\end{aligned}
\end{equation*}
where $z(t)$ is defined by \eqref{BA29}, and \eqref{BA30} now
looks
\begin{equation}\label{QL20}
\begin{aligned}
&x_j(\sigma_1\wedge t)=x_j(\tau_{l-1})+\left(1-\sum_{j\neq
j_0,j_1,\ldots,j_{l-2}}x_j(\tau_{l-1})\right)\times\\
&\times\left((\lambda_j(\tau_{l-1})-\mu_j)((\sigma_1\wedge
t)-\tau_{l-1})-\lambda_j(\tau_{l-1}) \int_0^{(\sigma_1\wedge
t)-\tau_{l-1}}z(s)\mbox{d}s\right),
\end{aligned}
\end{equation}
where $z(t)$ is defined by \eqref{BA32}, and $l=1,2,\ldots,\ell$.
\end{thm}

\section{Application: Performance analysis of symmetric large
closed client/server computer networks with unreliable client
stations} \label{Application} In this section we discuss
application of the above main theorem for symmetric large closed
client server computer networks.

\subsection{Formulation of the problem} We consider a network with
$k$ identical servers. We assume that all of the processes started
at zero, i.e. $A_{j,N}(0)=S_{j,N}(0)=Q_{j,N}(0)=0$, and the
following other assumptions and notation are used. The service
time of each unit in the server station is exponentially
distributed with parameter $\lambda$. There are $k$ client
stations in total, and each of client stations is a subject to
breakdown. A lifetime of each client station has the probability
distribution $G(x)$. In this case the moments of breakdown are
associated with change states of semi-Markov environment, and this
example is one of applications of the above theory.

The process $S_{j,N}(t)$ satisfies the condition
\begin{equation}\label{A1}
\mathbb{P}\left\{\lim_{N\to\infty}\frac{S_{j,N}(t)}{Nt}=\mu\right\}=1.
\end{equation}
(The equivalent form of condition \eqref{A1} is considered in
Section \ref{Definitions}. It is assumed here that
$\mu_1=\mu_2=\ldots=\mu_k\equiv\mu$.)

The relations between parameters $\lambda$, $\mu$ and $k$ are
assumed to be
\begin{equation}\label{A2}
\frac{\lambda}{k\mu}<1,
\end{equation}
and
\begin{equation}\label{A3}
\frac{\lambda}{\mu}>1.
\end{equation}
Condition \eqref{A2} means that all of the client stations are
initially non-bottleneck. Condition \eqref{A3} means that after
one of other breakdown all of the client stations become
bottleneck. Denote
\begin{equation*}
l_0=\max\left\{l: \frac{\lambda}{l\mu}>1\right\}
\end{equation*}
the maximum number of available client stations when the client
stations all are bottleneck. Then for all $l\leq l_0$ the rest $l$
client stations will be bottleneck as well.

Let $\alpha<1$ be a given positive number. We say that the network
is \textit{at risk} if the total number of units in client
stations exceeds the value $\alpha N$. Assuming that at the
initial time moment all of the units are in the server stations,
the aim of this section is to find a confidence interval [0,
$\theta$) such that with a given high probability $P$ (say
$P=0.95$) the network will not be at risk during that time
interval [0, $\theta$). The parameter $N$ is assumed to be large.
Therefore we study asymptotic solution of this problem as
$N\to\infty$, that is we study limiting (as $N\to\infty$)
normalized cumulative queue-length process in client stations.
This problem is hard in general. Here we consider a particular
case of $l_0=1$.

\subsection{Solution of the problem} The limiting normalized
queue-length process is denoted $q(t)$. At the initial time moment
$t=0$ there are $k$ available client stations. Let $\tau_1$,
$\tau_2$,\ldots,$\tau_k$ be the moments of their breakdown, $0\leq
\tau_1\leq\tau_2\leq\ldots\leq\tau_k$. The above moments of
breakdown are associated with the behavior of the time dependent
network, which can be considered as a network in semi-Markov
environment. Let us then apply Theorems \ref{thm1} and \ref{thm2}.
(Last Theorem \ref{thm2} is used consequently from one random
interval to another. Here we do not explain the details of this
application assuming that they are clear.)

The random time interval [0, $\tau_k$] is the lifetime of the
entire system. Therefore $q(t)$ is to be considered during the
aforementioned random interval [0, $\tau_k$]. We assume that
$l_0=1$. Therefore, according to Theorem \ref{thm1} and Theorem
\ref{thm2} we obtain that in the random interval [0,
$\tau_{k-1}$), $q(t)=0$, while in the random interval
[$\tau_{k-1}$, $\tau_k$) the equation for $q(t)$ is
\begin{equation}\label{A4}
q(t)=(\lambda-\mu)(t-\tau_{k-1})-\lambda\int_{0}^{t-\tau_{k-1}}r(s)\mbox{d}s,
\end{equation}
where $r(t)$ is given by
\begin{equation}\label{A5}
r(t)=\left(1-\frac{\mu}{\lambda}\right)(1-\mbox{e}^{-\lambda t}).
\end{equation}
In the last endpoint $\tau_k$ we set $q(\tau_k)=1.$

We have the following relationships:
\begin{equation}\label{A6}
\mathbb{P}\{q(t)=0\}=\sum_{i=2}^k\binom{k}{i}[1-G(t)]^i[G(t)]^{k-i},
\end{equation}
\begin{equation}\label{A7}
\begin{aligned}
&\mathbb{P}\{q(t)\leq\gamma<1\}\\&=[1-G(t)]\sum_{i=1}^{k-1}\binom{k-1}{i}[1-G(t-t_\gamma)]^i[G(t-t_\gamma)]^{k-i-1},
\end{aligned}
\end{equation}
where $t_\gamma$ is such the value of $t$ under which
\begin{equation}\label{A7+}
(\lambda-\mu)t-\lambda\int_0^tr(s)\mbox{d}s=\gamma.
\end{equation}

The value $t_\gamma$ is found from the relation
\begin{equation}\label{A8}
\frac{\int_0^\infty[1-G(t)]\sum_{i=1}^{k-1}\binom{k-1}{i}[1-G(t-t_\gamma)]^i[G(t-t_\gamma)]^{k-i-1}\mbox{d}t}
{\int_0^\infty\sum_{i=2}^{k}\binom{k}{i}[1-G(t-t_\gamma)]^{i}[G(t-t_\gamma)]^{k-i}\mbox{d}t}=P,
\end{equation}
and then from \eqref{A7+} one can find the corresponding value
$\gamma$.

If the value of $\gamma$ is not greater than $\alpha$, then the
value $\theta$ of the interval [$\tau_{k-1}$, $\theta$) is to be
taken $\theta=\tau_{k-1}+t_\gamma$. Otherwise, if $\gamma>\alpha$,
then the value $\theta$ is to be taken
$\theta=\tau_{k-1}+t_\alpha$.

In the particular case of $k=2$ we have the following results.
Relations \eqref{A6}, \eqref{A7} and \eqref{A8} reduces
correspondingly to
\begin{equation}\label{A9}
\mathbb{P}\{q(t)=0\}=[1-G(t)]^2,
\end{equation}
\begin{equation}\label{A10}
\begin{aligned}
&\mathbb{P}\{q(t)\leq\gamma<1\}=[1-G(t)][1-G(t-t_\gamma)],
\end{aligned}
\end{equation}
and
\begin{equation}\label{A11}
\frac{\int_0^\infty[1-G(t)][1-G(t-t_\gamma)]\mbox{d}t}
{\int_0^\infty[1-G(t-t_\gamma)]^2\mbox{d}t}=P.
\end{equation}

\subsection{Numerical calculation} We consider the following
example for $k=2$: $\lambda=4$, $\mu=3$, $\alpha=0.2$, $P=0.95$,
$G(x)=1-\mbox{e}^{-2x}$. From \eqref{A11} we have:
\begin{equation*}
\frac{\int_0^\infty\mbox{e}^{-2(t-t_\gamma)}\mbox{e}^{-2t}\mbox{d}t}{\int_0^\infty\mbox{e}^{-4(t-t_\gamma)}\mbox{d}t}
=\mbox{e}^{-2t_\gamma}=0.95.
\end{equation*}
Solution of the equation $\mbox{e}^{-2t_\gamma}=0.95$ yields
$t_\gamma=0.025647$. From \eqref{A7+} we obtain:
\begin{equation*}
\gamma=\int_0^{t_\gamma}\mbox{e}^{-4t}\mbox{d}t=\int_0^{0.025647}\mbox{e}^{-4t}\mbox{d}t=0.25-0.25\mbox{e}^{-0.102588}
\approx0.024375.
\end{equation*}
This value of $\gamma$ is less than $\alpha=0.2$, and therefore
this value $\gamma=0.024375$ is a required value for parameter,
which defines a desired confidence interval.

\section{Discussion of new problems and associated monotonicity
conditions for the networks in Markov
environment}\label{Discussion} Theorem \ref{thm2} looks very
complicated, and its further analysis is very difficult to make a
conclusion on the behavior of queue-length processes. For example,
it seems very difficult to obtain any numerical characteristics of
normalized queue-length processes analytically,
$\mathbb{E}x_j(\sigma_1\wedge t)$ for example. Therefore, the
numerical work should be based on simulation of semi-Markov
environment in order to obtain required performance
characteristics of the process. (By simulation of semi-Markov
environment we mean a multiple realizations for a Markov process
in order to calculate required numerical performance
characteristics of the process.)

For the purpose of performance analysis we also should restrict
the class of networks and processes describing the behavior of
queue-length processes in client stations. This restriction is
related to application of the results rather than development of
the theory. In many practical examples the quality characteristics
of networks are changed monotonically resulting in one or other
strategy of repair mentioned in Section \ref{Motivation}.

The aforementioned comparing of two different strategies for fixed
interval (0, $T$) requires application of Theorem \ref{thm2}, and
the problem can be solved without any additional assumption
requiring monotonicity. However, under general settings we cannot
answer to many significant questions. One of them is \textit{How
behave this criteria when the considered time interval is
changed?} For example, we have two strategies corresponding two
different initial conditions of semi-Markov environment, and
suppose we concluded that the first strategy is more profitable
than the second one for specific time interval (0, $T$).
\textit{Is this conclusion remains correct (or becomes not
correct) for another time interval (0, $T^{*}$)?} Another typical
question is as follows. Again, we have two strategies
corresponding two different initial conditions of semi-Markov
environment. Suppose we established that for an interval (0, $T$)
the both strategies are equivalent. Let $T^*$ be a new time
instant, and $T^*>T$. \textit{Which one of the strategies is now
more profitable in the new time interval (0, $T^*$), the first or
second one?} The same question can be asked under the opposite
inequality $T^*<T$.

These questions can be answered in the case when the class of the
processes studied numerically has a monotone stricture and is
described by the properties listed below. Then in certain cases
the behavior of queue-length processes in client stations and
consequently a conclusion about better strategy for other time
intervals can be established as well.

\smallskip
Assuming for simplicity that the environment is \textit{Markov},
then aforementioned properties are as follows.

\smallskip
(1) For any two positive integers $l\leq m$ assume that
$z_{l,m}\geq z_{m,l}$.

Recall that $z_{l,m}\triangle t$ + $o(\triangle t)$ ($l\neq m$)
are the transition probabilities  from the state $\mathcal{E}_l$
to the state $\mathcal{E}_m$ of a homogeneous Markov process for a
small time interval ($t$, $t+\triangle t$).

\smallskip
(2) $\lambda_j(\mathcal{E}_l)\leq \lambda_j(\mathcal{E}_m)$ for
all $j=1,2,\ldots,k$, and $l\leq m$.

\smallskip
Property (1) means that the Markov process $Z(t)$ is an increasing
process in the following sense: for two time moments
$\sigma_{l-1}$ and $\sigma_{l}$ we have
$Z(\sigma_{l-1})\leq_{st}Z(\sigma_l)$, which means that the state
of a Markov process in time $\sigma_{l-1}$ is not greater (in
stochastic sense) than that state in time $\sigma_{l}$ for any
integer positive $l$. The above property remains correct for any
$t_1\leq t_2$, i.e. $Z(t_1)\leq_{st}Z(t_2)$.

Property (1) also means that for two Markov processes $Z_1(t)$ and
$Z_2(t)$ having the same transition probabilities, but different
initial conditions satisfying $Z_1(0)\leq_{st} Z_2(0)$, we also
have $Z_1(t)\leq_{st} Z_2(t)$, $t\geq 0$. (For details of the
proof of these properties see e.g. Kalmykov \cite{Kalmykov 1962}.)

Consequently, from property (2) we have
$\lambda_j(\mathcal{E}(\sigma_{l-1}))\leq_{st}\lambda_j(\mathcal{E}(\sigma_{l}))$
for all $j=1,2,\ldots,k$ and any integer positive $l$. Moreover,
for all $j=1,2,\ldots,k$ and any $t_1\leq t_2$ we have
$\lambda_j(\mathcal{E}(t_1))\leq_{st}\lambda_j(\mathcal{E}(t_2))$.

Thus the rates $\lambda_j(\mathcal{E}(t))$ are increasing in time.
As a result, the queue-length processes in client stations
increase sharper than in the case of fixed $\lambda_j$ of ``usual"
network, and more extended problems mentioned in this section seem
can be solved as well. We however do not provide their solutions
in the present paper.

In the next section, numerical investigation for concrete client
server networks in a given Markov environment, satisfying the
above two properties is provided.

\section{Example of numerical study}\label{Numerical example}
In this section we do not intend to challenge a problem of
comparing two different strategies or finding an optimal strategy.
We only show (step-by-step) how to study the behavior of
queue-lengths in client stations numerically. However, the
detailed explanations of the given example can help to solve some
of the aforementioned problems of Section \ref{Motivation}. The
example models a Markov environment, i.e. in our example we
simulate exponentially distributed random variables describing the
state changes in the Markov environment. (In the given case by
simulation we mean one realization of the process in order to
study numerically a specific sample path of normalized
queue-length process.)

We consider the simplest case of Markov transition matrix of the
order 4
\begin{equation*}\label{Ex1}
P=
\left(\begin{matrix} 0 &1 &0 &0\\
0 &0 &1 &0\\
0 &0 &0 &1\\
0 &0 &0 &1
\end{matrix}
\right)
\end{equation*}
associated with the continuous Markov process $Z(t)$. In this
matrix $P_{l,l+1}=1$, $l=1,2,3$, and $P_{4,4}$=1. This means that
if the initial state of the process is
$\mathcal{E}(0)=\mathcal{E}_1$, then the next state is
$\mathcal{E}(\sigma_1)=\mathcal{E}_2$. Consequently,
$\mathcal{E}(\sigma_2)=\mathcal{E}_3$, and
$\mathcal{E}(\sigma_3)=\mathcal{E}_4$. Then
$\mathcal{E}(\sigma_l)=\mathcal{E}_4$ for all $l\geq 3$. Assume
also that $z_{l,l+1}=1$, $l=1,2,3$, so that
$\mathbb{E}(\sigma_l-\sigma_{l-1})$ =1, \ $l$=1,2,3.

For simplicity, the network contains only 2 client stations.
Assume that $\lambda_1(\mathcal{E}_1)$=1,
$\lambda_1(\mathcal{E}_2)$=2, and for $l$=3,4,
$\lambda_1(\mathcal{E}_l)$=$3$. We also assume that
$\lambda_1(\mathcal{E}_l)$ = $\lambda_2(\mathcal{E}_l)$ for
$l=1,2,3,4$. The values $\mu_1$ = $\mu_2$ = 2. Next, $\beta_1$ =
$\beta_2$ = 0.1, i.e. at the initial time moment each client
station contains 10\% of all units in the queue.

We set $T$=3, and study behavior of queue-length processes in
client station in the time interval (0, 3). By simulation we
obtained the following exponentially distributed random variables:
0.5488, 1.0892 and 1.8734. The sum of these 3 random variables is
greater than 3, so this quantity of exponentially distributed
random variables is enough for our experiment.

Note, that $\lambda_1(0)=\lambda_2(0)=0.8$. There three time
intervals: [0, 0.5488), [0.5488, 1.6380), [1.6380, 3).

For the time interval [0, 0.5488) the two client stations are
absolutely non-bottleneck because
$\frac{\lambda_j(\mathcal{E}_1)}{\mu_j}$=0.5, and we have the
following equations:
\begin{equation}\label{Ex2}
\begin{aligned}
&x_1(t)=x_2(t)=0.1-0.96t-0.8\int_0^tz(s)\mbox{d}s,\\
&z(t)=-1.5\left(1-\mbox{e}^{-1.6t}\right).
\end{aligned}
\end{equation}
Therefore, from \eqref{Ex2} we obtain:
\begin{equation}\label{Ex3}
x_1(t)=x_2(t)=-0.65+0.24t+0.75\mbox{e}^{-1.6t}.
\end{equation}
Substituting 0.5488 for $t$ in \eqref{Ex3} one can see that
\begin{equation*}\label{Ex4}
x_1(0.5488)=x_2(0.5488)\approx-0.2066.
\end{equation*}
The endpoints are negative, therefore we are to find such the
values $\tau_1$ and $\tau_2$ such that $x_1(\tau_1)$=0 and
$x_2(\tau_2)$=0. In our case $\tau_1$ = $\tau_2\approx0.117$. This
means that $q_1(t)=q_2(t)=0$ for all $t\geq0.117$ of the given
interval [0, 0.5488). Therefore, in the endpoint of this interval
$q_1(0.5488)=q_2(0.5488)=0$.

Consider now the time interval [0.5488, 1.6380). In point 0.5488
we now set $x_1(0.5488)=x_2(0.5488)$=0. Therefore,
$\lambda_1(0.5488)$ =$\lambda_2(0.5488)$ =2, and the both client
stations are bottleneck in [0.5488, 1.6380). Since
$\frac{\lambda_j(0.5488)}{\mu_j}$=1, $j=1,2$, then $x_1(t)$ and
$x_2(t)$ are equal to zero in this interval, and $q_1(t)$ =
$q_2(t)$ =0 in this interval as well.

We arrive at the last time interval [1.6380, 3). Similarly to the
above, we have $x_1(1.6380)=x_2(1.6380)$=0, and
$\lambda_1(1.6380)=\lambda_2(1.6380)$=3, and the both client
stations are bottleneck in [1.6380, 3). However, in the both
client stations we have $\frac{\lambda_j(1.6380)}{\mu_j}$=1.5,
$j=1,2$. Therefore, after a little algebra we have the following
equations:
\begin{equation*}\label{Ex5}
q_j(t)=\frac{1}{6}\left(1-\mbox{e}^{-6(t-1.6380)}\right), \ j=1,2,
\end{equation*}
for all $t$ from the interval [1.6380, 3).

\section{Concluding remarks}\label{Concluding remarks}
In the present paper we introduced a class of client/server
networks in order to study performance measures of real
client/server networks. Our analysis was based on the results of
earlier papers related to closed queueing networks with
bottleneck. However, for purpose of real applications, we
developed the earlier results and provided complete analysis of
standard bottleneck client/server networks. We then extended our
results for client/server networks in semi-Markov environment. The
results obtained in this paper are then used for analysis of
confidence intervals of client/server networks with failing client
stations. Numerical study given in this paper will help to clearly
understand solution for many related problems. The future work can
be related to application of the theoretical results of this paper
to concrete technological problems similar to those formulated in
Section \ref{Motivation}.

\section*{Acknowledgement} The research was supported by
the Australian Research Council, grant \#DP0771338.



\begin{thebibliography}{10}
\bibitem{Abramov 2000}\textsc{Abramov, V.M.} (2000). A large
closed queueing network with autonomous service and bottleneck.
\emph{Queueing Systems}, 35, 23-54.

\bibitem{Abramov 2001}\textsc{Abramov, V.M.} (2001). Some results
for large closed queueing networks with and without bottleneck:
Up- and down-crossings approach. \emph{Queueing Systems}, 38,
149-184.

\bibitem{Abramov 2004}\textsc{Abramov, V.M.} (2004). A large
closed queueing network containing two types of node and multiple
customers classes: One bottleneck station. \emph{Queueing
Systems}, 48, 45-73.

\bibitem{Abramov 2005}\textsc{Abramov, V.M.} (2005). The stability
of join-the-shortest-queue models with general input and output processes.
arXiv: math/PR 0505040.

\bibitem{Abramov 2006}\textsc{Abramov, V.M.} (2006). The effective bandwidth
problem revisited. arXiv: math/PR 0604182.

\bibitem{Abramov reliability 2007}\textsc{Abramov, V.M.} (2007).
Confidence intervals associated with performance analysis of
symmetric large closed client/server computer networks.
\emph{Reliability: Theory and Applications}, 2 (2), 35-42.

\bibitem{Abramov reliability 2007-2}\textsc{Abramov, V.M.} (2007).
Further analysis of confidence intervals for large client/server
computer networks, 
submitted.

\bibitem{Anulova and Liptser 1990}\textsc{Anulova, S.V. and
Liptser, R. Sh.} (1990). Diffusion approximation for processes
with normal reflection. \emph{Theory of Probability and its
Applications}, 35, 413-423.

\bibitem{Baccelli and Makovsky 1986} \textsc{Baccelli, F. and
Makovsky, A.M.} (1986). Stability and bounds for single-server
queue in a random enviromnent. \emph{Stochastic Models}, 2,
281-292.

\bibitem{Berger Bregman and Kogan 1999}\textsc{Berger, A., Bregman, L. and Kogan, Ya.} (1999)
Bottleneck analysis in multiclass closed queueing networks and its
application. \emph{Queueing Systems}, 31, 217-237.

\bibitem{Borovkov 1976}\textsc{Borovkov, A.A.} (1976).
\emph{Stochastic Processes in Queueing Theory}. Springer, Berlin.

\bibitem{Borovkov 1984}\textsc{Borovkov, A.A.} (1984). \emph{Asymptotic Methods in
Queueing Theory}. John Wiley, New York.

\bibitem{Boxma and Kurkova 2000}\textsc{Boxma, O.J. and Kurkova,
I.A.} (2000). The M/M/1 queue in heavy-tailed random environment.
\emph{Statistica Neerlandica}, 54, 221-236.

\bibitem{Dellacherie 1972}\textsc{Dellacherie, C.} (1972).
\emph{Capacit\'es et Processus Stochastiques.} Springer-Verlag,
Berlin.

\bibitem{D'Auria 2007}\textsc{D'Auria, B.} M/M/$\infty$ queues in
quasi-Markovian random environment. arXiv: math/PR 0701842.

\bibitem{Fricker 1986}\textsc{Fricker, C.} (1986).  Etude d'une file GI/G/1
\'a service autonome (avec vacances du serveur).  \emph{Advances
in Applied Probability}, 18, 283-286.

\bibitem{Fricker 1987}\textsc{Fricker, C.} (1987). Note sur un modele de file GI/G/1
\'a service autonom\'e (avec vacances du serveur).  \emph{Advances
in Applied Probability}, 19, 289-291.

\bibitem{Gelenbe and Iasnogorodski 1979}\textsc{Gelenbe, E. and
Iasnogorodski, R.} (1980). A queue with server of walking type
(autonomous service). \emph{Ann. Inst. H. Poincare}, 16, 63-73.

\bibitem{Helm and Waldmann 1984}\textsc{Helm, W.E. and Waldmann,
K.-H.} (1984). Optimal control of arrivals to multiserver queues
in a random environment. \emph{Journal of Applied Probability},
21, 602-615.



\bibitem{Kalmykov 1962}\textsc{Kalmykov, G.I.} (1962). On the
partial ordering of one-dimensional Markov processes. \emph{Theory
of Probability and its Applications}, 7, 456-459.

\bibitem{Kogan 1992}\textsc{Kogan, Ya.} (1992). Another approach to
asymptotic expansions for large closed queueing networks.
\emph{Operations Research Letters}, 11, 317-321.



\bibitem{Kogan and Liptser 1993}\textsc{Kogan, Ya. and Liptser,
R.Sh.} (1993). Limit non-stationary behavior of large closed
queueing networks with bottlenecks. \emph{Queueing Systems}, 14,
33-55.

\bibitem{Kogan Liptser and Smorodinskii 1986}\textsc{Kogan, Ya.,
Liptser, R.Sh. and Smorodinskii, A.V.} (1986). Gaussian diffusion
approximation of a closed Markov model of computer networks.
\emph{Problems of Information Transmission}, 22, 38-51.

\bibitem{Krichagina Liptser and Puhalskii 1988}\textsc{Krichagina,
E.V., Liptser, R.Sh. and Puhalskii, A.A.} (1988). Diffusion
approximation for a system that an arrival stream depends on queue
and with arbitrary service. \emph{Theory of Probability and its
Applications}, 33, 114-124.

\bibitem{Krichagina and Puhalskii 1997}\textsc{Krichagina, E.V. and
Puhalskii, A.A.} (1997). A heavy-traffic analysis of closed
queueing system with $GI/\infty$ server. \emph{Queueing Systems},
25, 235-280.

\bibitem{Krieger et al 2005}\textsc{Krieger, U., Klimenok, V.I., Kazmirsky, A.V., Breuer, L.
and Dudin, A.N.} (2005). A BMAP/PH/1 queue with feedback operating
in a random environment. \emph{Mathematical and Computer
Modelling}, 41, 867-882.

\bibitem{Liptser and Shiryayev 1989}\textsc{Liptser, R.Sh. and
Shiryayev, A.N.} (1989). \emph{Theory of Martingales}. Kluwer,
Dordrecht.

\bibitem{McKenna and Mitra 1982}\textsc{McKenna, J. and Mitra, D.}
(1982). Integral representation and asymptotic expansions for
closed Markovian queueing networks. Normal usage. \emph{Bell
System Technical Journal}, 61, 661-683.

\bibitem{O'Cinneide and Purdue 1986}\textsc{O'Cinneide, C. and
Purdue, P.} (1986). The M/M/$\infty$ queue in a random
environment. \emph{Journal of Applied Probability}, 23, 175-184.

\bibitem{Pittel 1979}\textsc{Pittel, B.} (1979). Closed
exponential networks of queues with saturation: The Jackson type
stationary distribution and its asymptotic analysis.
\emph{Mathematics of Operations Research}, 6, 357-378.

\bibitem{Purdue 1974}\textsc{Purdue, P.R.} (1974). The M/M/1 queue in a random environment.
\emph{Operations Research}, 22, 562-569.

\bibitem{Ramanan 2006}\textsc{Ramanan, K.} (2006). Reflected
diffusions defined via extended Skorokhod map. \emph{Electronic
Journal of Probability}, 11, 934-992.

\bibitem{Skorokhod 1961}\textsc{Skorokhod, A.V.} (1961).
Stochastic equations for diffusion processes in a bouded region.
\emph{Theory of Probability and its Applications}, 6, 264-274.

\bibitem{Tanaka 1979}\textsc{Tanaka, H.} (1979). Stochastic
differential equations with reflected boundary conditions in
convex regions. \emph{Hiroshima Mathematical Journal}, 9, 163-177.

\bibitem{Whitt 1984}\textsc{Whitt, W.} (1984). Open and closed
models for networks of queues. \emph{AT\&T Bell Laboratory
Technical Journal}, 63, 1911-1979.


\end{thebibliography}
\end{document}